\documentclass[12pt]{amsart}

\usepackage{a4}
\usepackage{amssymb,amsfonts,amsmath}

\usepackage{color}

\makeatletter
\def\LaTeX{\leavevmode L\raise.42ex
\hbox{\kern-.3em\size{\sf@size}{0pt} \selectfont
A}\kern-.15em\TeX} \makeatother

\newcommand{\BibTeX}{{\rm B\kern-.05em{\sc i\kern-.025emb}
\kern-.08em\TeX}}

\theoremstyle{definition}

\makeatletter\label{e:dispaa}
\label{e:dispau}
\label{e:dispav}
\label{e:dispaw}
\label{e:dispax}
\def\theequation{\thesection.\@arabic\c@equation}
\makeatother\makeatletter

\newcommand{\vs}{\vspace{.3 cm}}
\newcommand{\vsd}{\vspace{.2 cm}}
\newcommand{\pn}{\par\noindent}
\newcommand{\lra}{\longrightarrow}
\newcommand{\ra}{\rightarrow}

\newcommand{\benu}{\begin{enumerate}}
\newcommand{\enu}{\end{enumerate}}

\makeatother
\newcommand{\A}{{\mathbb{A}}}
\newcommand{\D}{{\mathbb{D}}}
\newcommand{\E}{{\mathbb{E}}}

\newcommand{\N}{{\mathbb{N}}}

\newcommand{\T}{{\mathbb{T}}}

\newcommand{\bema}{\left ( \begin{array}}
\newcommand{\ema}{\end{array} \right )}

\begin{document}
\title[Representation dimension]%
{The representation dimension of a class of tame algebras}

\author[Assem]{Ibrahim Assem}
\address{D\'epartement de math\'ematiques, Universit\'e de Sherbrooke,
Sherbrooke, Qu\'ebec, Canada, J1K 2R1.}
\email{ibrahim.assem@usherbrooke.ca}

\author[Coelho]{Fl\'avio U. Coelho}
\address{Departamento de Matem\'atica-IME, Universidade de S\~ao Paulo,
CP 66281, S\~ao Paulo, SP, 05315-970, Brazil}
\email{fucoelho@ime.usp.br}

\author[Trepode]
{Sonia Trepode}
\address{Depto de Matem\'atica, FCEyN, Universidad Nacional de Mar del Plata, 7600,
Mar del Plata, Argentina} \email{strepode@mdp.edu.ar}

\thanks{ This paper was finished when the second  author was
visiting the Universit\'e de Sherbrooke under the support of NSERC, to whom he is very 
grateful.  The first author gratefully acknowledges partial support from the NSERC 
of Canada and the Universit\'e de Sherbrooke. The second author thanks the CNPq. 
The third author is a researcher of CONICET. The authors thank Grzegorz Bobinski 
for pointing out an inacurrancy in a previous version.}

\maketitle

\begin{abstract}
We prove that, if $A$  is a strongly simply 
connected algebra of polynomial growth then $A$ is torsionless-finite. In particular, 
its representation  dimension is at most three. 
\end{abstract}
\vspace{.3 cm}

Among the most useful algebraic invariants are the homological dimensions, which are meant to 
measure how much an algebra or a module deviates from a situation considered to be ideal. Introduced 
by Auslander in the early seventies \cite{Au}, the representation dimension of an Artin algebra 
was long left aside from the mainstream of the theory, until a marked renewal of interest about 
ten years ago. It measures the least global dimension of all endomorphism rings  of those finitely 
generated modules which are both generators and cogenerators of the module category. Part of the 
reason for this new interest comes from the fact that Igusa and Todorov have shown that, if the 
representation dimension of an Artin algebra is at most three, then its finitistic dimension is 
finite \cite{IT}. Also, Iyama has proved that, for any Artin algebra $A$, the representation 
dimension rep.dim.$A$ is finite \cite{Iy}. Since then, there have been several attempts to 
understand this invariant and to calculate it for classes of algebras, see for instance \cite{APT,CHU,CP,Op}
or the survey \cite{Ri3}. It was shown by Auslander that an Artin algebra is representation-finite 
if and only if its representation dimension is at most two \cite{Au}. Since Auslander's expectation 
was that this dimension would measure how far an algebra is from being representation-finite, 
it is natural to ask whether tame algebras have representation dimension at most three. The answer 
to this question is known to be positive for some classes of tame algebras, such as special 
biserial algebras \cite{EHIS} and domestic self-injective algebras socle equivalent to a weakly 
symmetric algebra of euclidean type \cite{BHS}.

The objective of our paper is to prove that the representation dimension of a strongly simply 
connected algebra of polynomial growth (over an algebraically closed field) is at most three. 
These algebras form a nice class which has been extensively studied by Skowro\'nski, de la Pe\~na 
and others, see, for instance, the survey \cite{AS4}. In particular, it is shown in \cite{Sk2} that 
strongly simply connected algebras of polynomial growth are multicoil algebras. We recall also 
that, as pointed out in \cite{Ri3}, an algebra which is torsionless-finite (that is, such that 
any indecomposable projective module has only finitely many isomorphism classes of indecomposable 
submodules) has its representation dimension at most three. Our main result is the following theorem. 
\vs\\
{\sc Theorem A.} {\it Let $A$ be a strongly simply connected algebra of polynomial growth. Then $A$ is torsionless-finite. In particular, rep.dim$.A \leq 3$.}
\vs

In the course of the proof, we found that the following result, of independent interest, was necessary. 
\vs\\
{\sc Theorem B.} {\it Let $C$ be a tame concealed algebra of type distinct of ${\widetilde{\A}}$, and $(H_{\lambda})_\lambda$ be an infinite family 
of pairwise non-isomorphic simple homogeneous $C$-modules. Suppose $M$ is a $C$-module such that 
$H_{\lambda} \subset M$ for all $\lambda$. Then $C[M]$ is wild. } 
\vs

The paper is organized as follows. After a short preliminary section, in which we fix the notations 
and recall facts about multicoil algebras, we prove our Theorem(B) in section 2 and Theorem(A) in section 3.

\section{Multicoil algebras}

\subsection{Notation.} 
In this paper, $k$ denotes a fixed algebraically closed field. By an algebra $A$ is meant a basic, 
connected, associative finite dimensional $k$-algebra with an identity. Thus, there exists a connected
bound quiver $(Q_A, I)$ and an isomorphism $A \simeq kQ_A/I$. Equivalently, $A$ may be viewed as a 
$k$-category, of which the object class $A_0$ is the set of points of $Q_A$ and the set of morphisms 
from $x$ to $y$ is the quotient of the $k$-vector space $kQ_A(x,y)$ of linear combinations of paths 
in $Q_A$ from $x$ to $y$ by $I(x,y) = I\cap kQ_A(x,y)$, see \cite{BG}. A full subcategory $C$ of $A$ is 
{\it convex} (in $A$) if, for any path $x_0 \ra x_1 \ra \cdots \ra x_t$ in $A$ with $x_0, x_t \in C_0$, 
we have $x_i \in C_0$ for all $i$. The algebra $A$ is {\it triangular} if $Q_A$ is acyclic. 

By $A$-module is meant a finitely generated right $A$-module. We denote by mod$A$ the category of $A$-modules
and by ind$A$ a full subcategory consisting of a complete set of representatives of the isomorphism 
classes of indecomposable $A$-modules. We recall that, if $A \simeq kQ/I$, then an $A$-module $M$ is 
identified to a corresponding representation $(M(x)_{x\in Q_0}, M(\alpha)_{\alpha\in Q_1})$ of the 
bound quiver $(Q,I)$, see \cite{ASS}. For a point $x\in Q_0$, we denote by $P_x$ (or $I_x$, or $S_x$) the 
indecomposable projective (or injective, or simple, respectively) $A$-module corresponding to $x$. The 
{\it support} of an $A$-module $M$ is the full subcategory Supp$M$ of $A$ with objects those $x\in A_0$ 
such that $M(x) \neq 0$. For a full subcategory ${\mathcal C}$ of mod$A$, we denote by add${\mathcal C}$ 
the additive full subcategory with objects the direct sums of direct summands of objects in ${\mathcal C}$. 
If ${\mathcal C}$ contains a single module $M$, we write add${\mathcal C} =$ add$M$. For two full 
subcategories ${\mathcal C}, {\mathcal C}'$ of ind$A$, the notation Hom$_A({\mathcal C}, {\mathcal C}') =0$ 
means that Hom$_A(M, M')= 0$ for all $M\in {\mathcal C}, M\in  {\mathcal C}'$. We then denote by 
${\mathcal C}'\vee {\mathcal C}$ the full subcategory of ind$A$ having as objects those of 
${\mathcal C}_0' \cup {\mathcal C}_0$. 

A {\it path} in ind$A$ from $M$ to $N$ (sometimes denoted as $M \leadsto N$) is a sequence of non-zero 
morphisms
$$  M = M_0 \lra M_1 \lra \cdots \lra M_t = N \leqno(*) $$
with all $M_i$ in ind$A$. A path $(*)$ is a {\it cycle} if $M=N$ and at least one of the morphisms is not 
an isomorphism. An indecomposable module is {\it directed} if it lies on no cycle. 

We use freely properties of the Auslander-Reiten translations $\tau_A =$ DTr and $\tau^{-1}_A = $ TrD 
and the Auslander-Reiten quiver $\Gamma($mod$A$) of $A$ for which we refer to \cite{ASS,Ri2}. We identify points 
in $\Gamma($mod$A$) with the corresponding indecomposable $A$-modules, and components with the 
corresponding full subcategories of ind$A$. A component $\Gamma$ is {\it standard} if the category 
$\Gamma$ is equivalent to the mesh category $k(\Gamma)$. For tubes, tubular extensions and 
coextensions, we refer to \cite{Ri2} and for tame algebras, we refer to \cite{Sk1}. 

\subsection{One-point extensions.} 
The {\it one-point extension} of an algebra $A$ by an $A$-module $M$ is the matrix algebra 
$$ A[M] \ = \ \bema{cc} A & 0 \\ M & k \ema $$
with the usual addition and multiplication of matrices. The quiver of $A[M]$ contains $Q_A$ as a full 
convex subquiver and there is an additional (extension) point which is a source. The $A[M]$-modules 
are identified with triples $(V, X, \varphi)$, where $V$ is a $k$-vector space, $X$ an $A$-module and 
$\varphi \colon V \lra $ Hom$_A(M,X)$ is a $k$-linear map. An $A[M]$-linear map 
$(V,X, \varphi) \lra (V', X', \varphi ')$ is a pair $(f,g)$ where $f \colon V \lra V'$ is $k$-linear 
and $g \colon X \lra X'$ is $A$-linear such that $\varphi'f = $ Hom$_A(M,g)\varphi$. The dual notion 
is that of {\it one-point coextension}. 

A {\it vector space category} \cite{Ri1,Ri2} ${\mathcal K}$ is a $k$-category together with a faithful 
$k$-linear functor $|\cdot | \colon {\mathcal K} \lra $ mod$k$. The {\it subspace category} 
${\mathcal U}({\mathcal K})$ of ${\mathcal K}$ has as objects the triples $(V, X, \varphi)$ where $V$ is a 
$k$-vector space, $X$ an object in ${\mathcal K}$ and $\varphi \colon V \lra |X|$ is a $k$-linear 
monomorphism. A morphism $(V,X, \varphi) \lra (V', X', \varphi ')$ is a pair $(f,g)$ where 
$f \colon V \lra V'$ is $k$-linear and $g \colon X \lra X'$ is a morphism in ${\mathcal K}$ 
such that $\varphi'f = |g|\varphi$. 

If $A$ is an algebra and $M$ an $A$-module, one considers the vector space category Hom$_A(M,$ mod$A$) 
whose objects are of the form Hom$_A(M,X)$ with $X$ an $A$-module, and morphisms are of the form 
$$ \mbox{Hom}_A(M,f) \colon \mbox{ Hom}_A(M,X) \lra \mbox{ Hom}_A(M',X'), $$ 
where $f \colon X \lra X'$ is $A$-linear. Then
$|$Hom$_A(M,X)|$ is just the underlying\break 
$k$-vector space of Hom$_A(M,X)$. It is shown in \cite{Ri1} that 
${\mathcal U}($Hom$_A(M,$mod$A$)) is equivalent to the full subcategory of mod$A[M]$ consisting of 
the triples $(V,X, \varphi)$ without non-zero direct summands of the form $(k,0,0)$ or $(0, Y,0)$ 
where\break 
Hom$_A(M,Y) = 0$. We need essentially the following lemma. 
\vs\\
{\sc Lemma.} {\it Let $A$ be a tame algebra, and $M$ be an $A$-module. If $L$ is a submodule of 
$M$ such that $A[L]$ is wild, then $A[M]$ is wild.} 
\begin{proof} 
Since $A$ is tame, while $A[L]$ is wild, then the vector space category\break 
${\mathcal U}($Hom$_A(L,$mod$A$))
is wild. Now the inclusion $L \hookrightarrow M$ induces an epifunctor 
$$ {\mathcal U}(\mbox{Hom}_A(M,\mbox{mod}A)) \lra {\mathcal U}(\mbox{Hom}_A(L,\mbox{mod}A)).  $$
Since the latter is wild, then so is ${\mathcal U}(\mbox{Hom}_A(M,\mbox{mod}A))$. Therefore, 
$A]M]$ is wild.
\end{proof} 

\subsection{Coils.} 
Let $A$ be an algebra, and $\Gamma$ a standard component of $\Gamma($mod$A$). Given a module 
$X \in \Gamma$, called the {\it pivot}, the {\it support} Supp$(X, -)|_\Gamma$ of the functor 
Hom$_A(X, -)|_\Gamma$ is defined as follows. Let ${\mathcal H}_X$ be the full subcategory of 
ind$A$ consisting of the $Y \in \Gamma$ such that Hom$_A(X,Y) \neq 0$, and ${\mathcal I}_X$ 
be the ideal of ${\mathcal H}_X$ consisting of the morphisms $f \colon Y \lra Y'$ (with 
$Y, Y'$ in ${\mathcal H}_X$) such that Hom$_A(X, f) \neq 0$. Then 
Supp$(X, -)|_\Gamma = {\mathcal H}_X/ {\mathcal I}_X$. We define three admissible operations:
\benu
\item[(ad1)] Assume Supp$(X,-)|_\Gamma$ consists of an infinite sectional path starting at $X$
$$ X = X_0 \lra X_1 \lra X_2 \lra \cdots $$
Let $t\geq1$, $D = \T_t (k)$ be the full $t\times t$ lower triangular matrix algebra, and $Y$ 
be the unique indecomposable projective-injective $D$-module. The {\it modified algebra} of 
$A$ is $A' = (A\times D)[X \oplus Y]$. If $t=0$, it is simply $A' = A[X]$. 
\item[(ad2)] Assume Supp$(X,-)|_\Gamma$ consists of two sectional paths starting at $X$, the first
infinite and the second finite with at least one arrow
$$ Y_t \longleftarrow \cdots \longleftarrow Y_2 \longleftarrow Y_1 \longleftarrow X = X_0 
\lra X_1 \lra X_2 \lra \cdots $$
with $t\geq 1$. The {\it modified algebra} is $A' = A[X]$.
\item[(ad3)] Assume Supp$(X,-)|_\Gamma$ consists of two parallel sectional paths, the first 
infinite and starting at $X$, and the second finite with at least one arrow
\begin{picture}(60,19)
\put(21,0){$X = X_0$}
\put(37,0){$\lra$} 
\put(45,0){$X_1$}
\put(52,0){$\lra$}
\put(63,0){$\cdots$}
\put(71,0){$\lra$}
\put(79,0){$X_{t-1}$}
\put(87,0){$\lra$}
\put(95,0){$X_{t}$}
\put(102,0){$\lra$}
\put(112,0){$\cdots$}
\put(32,6){$\uparrow$}
\put(46,6){$\uparrow$}
\put(80,6){$\uparrow$}
\put(31,11){$Y_1$}
\put(37,11){$\lra$}
\put(45,11){$Y_2$}
\put(52,11){$\lra$}
\put(63,11){$\cdots$}
\put(71,11){$\lra$}
\put(79,11){$Y_t$}

\end{picture}
\vspace{.5 cm}\\
where $t\geq 2$. In particular, $X_{t-1}$ is injective. The {\it modified algebra} 
is $A' = A[X]$. 
\enu

In each case, $t$ is the {\it parameter} of the operation, and the component $\Gamma'$ of $\Gamma($mod$A'$) 
containing $X$ is the {\it modified component}. We also consider the duals of the above operations, 
denoted by (ad1*), (ad2*) and (ad3*), respectively. These six operations are the {\it admissible operations}. 

An Auslander-Reiten component $\Gamma$ is a {\it coil} if there exists a sequence\break 
$\Gamma_0, \Gamma_1, \cdots, \Gamma_m = \Gamma$ where $\Gamma_0$ is a stable tube and, for each $i$, 
$\Gamma_{i+1}$ is obtained from $\Gamma_i$ by an admissible operation.

Let $A$ be an algebra $A$. A family ${\mathcal R} = ( {\mathcal R}_{\lambda})_{\lambda \in \Lambda}$ of 
components of $\Gamma($mod$A$) is {\it weakly separating} if the indecomposable $A$-modules not in 
${\mathcal R}$ split into two classes ${\mathcal P}$ and ${\mathcal Q}$ such that 
\benu
\item[(1)] The ${\mathcal R}_\lambda$ are standard and pairwise orthogonal. 
\item[(2)] Hom$_A({\mathcal Q}, {\mathcal P}) $ = Hom$_A({\mathcal Q}, {\mathcal R}) $ = 
Hom$_A ({\mathcal R},{\mathcal P}) = 0$.
\item[(3)] Any morphism from ${\mathcal P}$ to ${\mathcal Q}$ factors through add${\mathcal R}$.
\enu

If ${\mathcal R}$ is a weakly separating family in $\Gamma$(mod$A$) consisting of stable tubes, then 
an algebra $B$ is a coil enlargement of $A$ using modules from ${\mathcal R}$ if there exists a finite sequence
of algebras $A = A_0, A_1, \cdots , A_m = B$ such that, for each $i$, $A_{i+1}$ is obtained from $A_i$ 
by an admissible operation with pivot either on a stable tube of ${\mathcal R}$ or on a coil of 
$\Gamma($mod$A_i)$ obtained from a stable tube of ${\mathcal R}$ by means of the admissible operations 
done so far. The {\it coil type} $c_B = ( c^-_B, c^+_B)$ of $B$ is a pair of functors $\Lambda \lra \N$ 
defined by induction on $i$, for each $\lambda \in \Lambda$, as follows.
\benu
\item[(a)] $c_A = (c^-_0, c^+_0)$ is such that $c^-_0(\lambda) = c^+_0(\lambda)$ is the rank of the 
stable tube ${\mathcal R}_\lambda$.
\item[(b)] If $c_{A_{i-1}} = (c^-_{i-1}, c^+_{i-1})$ is known and $t_i$ is the parameter of the 
operation from $A_{i-1}$ to $A_i$, then $c_{A_i} = (c^-_i, c^+_i)$ is defined by: 
$$ c^-_i(\lambda) \ = \ \left\{ \begin{array}{ll} c^-_{i-1}(\lambda) + t_i + 1 & 
\mbox{if the operations is (ad1*), (ad2*), (ad3*)}\\
& \mbox{with pivot in the coil of } \Gamma(\mbox{mod}A_{i-1})\\ 
&  \mbox{arising from } {\mathcal R}_\lambda .\\
c^-_{i-1}(\lambda) & \mbox{otherwise} \end{array} \right. $$
and
$$ c^+_i(\lambda) \ = \ \left\{ \begin{array}{ll} c^+_{i-1}(\lambda) + t_i + 1 & 
\mbox{if the operations is (ad1), (ad2), (ad3)}\\
& \mbox{with pivot in the coil of } \Gamma(\mbox{mod}A_{i-1})\\ 
&  \mbox{arising from } {\mathcal R}_\lambda .\\
c^+_{i-1}(\lambda) & \mbox{otherwise} \end{array} \right. $$
\enu

If all but at most finitely many values of each of $c^-_B, c^+_B$ equal 1, we replace each 
sequence by a finite sequence containing at least two terms, and including all those exceeding 1. 
\vsd\\
{\sc Theorem.} \cite{AST} {\it Let $A$ be an algebra having a weakly separating family of stable
tubes ${\mathcal R}$, and $B$ be a coil enlargement of $A$ using modules from ${\mathcal R}$. Then
\benu
\item[(a)] There is a unique maximal branch coextension $B^-$ (or extension $B^+$) of $A$ which is 
a full convex subcategory of $B$, and having $c^-_B$ as coextension type (or $c^+_B$ as extension 
type, respectively). 
\item[(b)] ind $B = {\mathcal P}' \vee {\mathcal R}' \vee {\mathcal Q}'$, where ${\mathcal R}'$ is a 
weakly separating family of ind$B$ obtained from ${\mathcal R}$ by the sequence of admissible operations, 
and separating ${\mathcal P}'$ from ${\mathcal Q}'$. Moreover, ${\mathcal P}'$ consists of $B^-$-modules, 
while ${\mathcal Q}'$ consists of $B^+$-modules.
\enu }

\subsection{Tame coil enlargements.} 
Let $B$ be a coil enlargement of a tame concealed algebra. Its {\it coil type} $c_B = (c^-_B, c_B^+)$ is 
tame if each of the sequences $c_B^-$, $c_B^+$ is one of the following: $(p,q)$ with $1 \leq p\leq q$, 
$(2,2,r)$ with $r \geq 2$, (2,3,3), (2,3,4), (2,3,5), (2,4,4), (2,3,6) or(2,2,2,2). We have the following 
result.
\vsd\\
{\sc Corollary.} \cite{AST}(4.3) {\it Let $B$ a coil enlargement of a tame concealed algebra. The 
following conditions are equivalent:
\benu
\item[(a)] $B$ is tame. 
\item[(b)] $B$ is of polynomial growth. 
\item[(c)] $c_B$ is tame.
\item[(d)] $B^-$ and $B^+$ are tame.
\enu
Moreover, $B$ is domestic if and only if $B^-$ and $B^+$ are tilted of euclidean type. }

\subsection{Multicoil algebras.} 
An Auslander-Reiten component $\Gamma$ is a {\it multicoil} if it contains a full translation 
subquiver $\Gamma'$ which is a disjoint union of coils such that no point in $\Gamma \setminus \Gamma'$ 
belongs to a cyclical path. 

An algebra $A$ is a {\it multicoil algebra} if, for any cycle $M_0 \ra M_1 \ra \cdots \ra M_t= M_0$
in ind$A$, all the $M_i$ lie in one standard coil of a multicoil of $\Gamma($mod$A$). The first part of 
the following theorem  is \cite{AS2} (4.6), the second part is \cite{Sk2}(4.1).
\vsd\\
{\sc Theorem.} \cite{AS2,Sk2} {\it Let $A$ be a multicoil algebra, then $A$ is of polynomial growth. If $A$ is strongly 
simply connected, then the converse also holds. }
\vs

If $A$ is a multicoil algebra, then it is triangular (hence of finite global dimension) \cite{AS3}(3.5). 
Also, any full convex subcategory of $A$ is a multicoil algebra \cite{AS3}(5.6). 

\subsection{Supports.} Supports of indecomposable modules over multicoil algebras are characterised 
in the following lemma. 
\vsd\\
{\sc Lemma.} Let $M$ be an indecomposable module over a multicoil algebra $A$, and let $B =$ Supp$M$. 
Then: 
\benu
\item[(a)] If $M$ is directed, then $B$ is tame and tilted. 
\item[(b)] If $M$ is not directed, then $B$ is a tame coil enlargement of a tame concealed algebra.
\enu
\begin{proof} Since $B$ is a full subcategory of $A$, which is tame, then $B$ is also tame. Then (a) 
follows from \cite{ASS}(IX.2.8, p. 366) and (b) follows from \cite{AS3}(5.9).
\end{proof}

Note that, if $L$ is an $A$-submodule of $M$, and $x \in A_0$ is such that $L(x) \neq 0$, 
then $M(x) \neq 0$. Hence, $L$ is also a $B$-submodule of $M$. 

\section{One-point extension of tame concealed algebras} 

\subsection{} 
In this section, we prove that, if $C$ is a tame concealed algebra and $M$ is a $C$-module containing an 
infinity of pairwise non-isomorphic simple homogeneous submodules, then $C[M]$ is wild. We start
with reduction lemmata. 
\vsd\\
{\sc Lemma.} {\it Let $C$ be a tame concealed algebra, let $M$ be a preinjective $C$-module, 
and $H$ a simple homogeneous submodule of $M$. Then:
\benu
\item[(a)] $M$ is sincere.
\item[(b)] If $C[M]$ is not wild, then dim$_k $(top$M) \leq 2$. Then, either $M$ is 
indecomposable, or $M \simeq M_1 \oplus M_2$, with $M_1, M_2$ indecomposables with simple top.
\enu } 
\begin{proof}
(a) This follows from the sincerity of $H$. 

(b) If dim$_k($top$M) =d  >2$, then $C[M]$ contains a wild full subcategory of the form\\
\begin{picture}(60,10)
\multiput(61,0)(22,0){2}{\circle*{1}}

\multiput(80,3)(0,-6){2}{\vector(-1,0){16}}
\multiput(72,2)(0,-1){5}{\circle*{.1}}

\put(70,5){$\alpha_1$}
\put(70,-8){$\alpha_d$}
\end{picture}
\vspace{1 cm}\\
a contradiction. The second statement follows from Nakayama's lemma, and (1.2). 
\end{proof}

\subsection{} 
Let $A = kQ/I$ be a bound quiver algebra. It is well-known (see, for instance, 
\cite{ASS}(III.2.2, p.77)) that, if $M$ is an $A$-module, then (top$M)(x) = M(x)$ is $x$ is a source,
and (top$M)(y) = \sum \{ $Coker$(M(\alpha)) \colon \alpha \colon x \lra y \}$, if $y$ is not a 
source.

If $M$ satisfies the hypothesis of (2.1)(b), and $x$ is a source, we then have dim$_kM(x) \leq 2$. 
\vsd\\
{\sc Corollary.}{\it Let $C$ be a tame concealed algebra, let $M$ be a preinjective $C$-module, 
and $H$ a simple homogeneous submodule of $M$. If $C[M]$ is not wild, then:
\benu
\item[(a)] $C$ has at most two sources.
\item[(b)] If $C$ has just one source $s$ and dim$_k H(s) = d \geq 2$, then $d=2$ 
and top$H$ = top$M$.
\item[(c)] If $C$ has two sources $s_1, s_2$ then dim$_kH(s_1) = $ dim$_k H(s_2) = 1$ and 
top$H = $ top$M$. 
\enu } 
\begin{proof} (a) If $C$ has at least three sources $s_1, s_2, s_3$, the remark above 
implies that top$M$ has $S_{s_1} \oplus S_{s_2} \oplus S_{s_3}$ as a summand, and 
dim$_k($top$M) \geq 3$, a contradiction. 

(b) Since $H(s) \subset M(s)$, which is at most two-dimensionsl, then $d=2$ and top$H = $top$M = S^2$.

(c) Clearly, $S_{s_1} \oplus S_{s_2}$ is a summand of both top$H$ and top$M$. Since 
dim$_k($top$M) \leq 2$, then top$M = S_{s_1} \oplus S_{s_2}$. We claim that top$H =  S_{s_1} \oplus S_{s_2}$.
If this is not the case, then there exists $a \notin \{ s_1, s_2 \}$ such that $S_a$ is a 
summand of top$H$. Hence (top$H)(a) \neq 0$. By the remark above, there exists an arrow 
$ \alpha \colon b \lra a$ in $C$ such that $H(\alpha)$ is not surjective. On the other hand, since $S_a$ is
not a summand of top$M$, then $M(\alpha)$ is surjective. Now there is a path from one of the 
$s_i$ (with $i \in \{ 1,2 \}$) to $a$ passing through $\alpha$, which we may assume, without loss of 
generality, to be of minimal length: 
$$s_i = b_t \lra b_{t-1} \lra \cdots \lra b_1 = b \stackrel{\alpha}{\lra} a \leqno(*)$$
We claim that dim$_kM(b) \geq 2$. Indeed, if this in not the case, then the inclusion 
$j \colon H \hookrightarrow M$ yields a commutative square\\
\begin{picture}(60,16)
\multiput(60,8)(0,-16){2}{\vector(1,0){24}}
\multiput(51.5,6)(38.5,0){2}{\vector(0,-1){10}}

\put(48,7){$H(b)$}   \put(87,7){$H(a)$}
\put(48,-9){$M(b)$}   \put(87,-9){$M(a)$}
\put(68,10){$H(\alpha)$}   \put(68,-6){$M(\alpha)$}
\put(46,-1){$j_b$}   \put(92,-1){$j_a$}
\end{picture}
\vspace{1.4 cm}\\
and $M(b) \simeq H(b) \simeq k$ gives that $j_b$ is an isomorphism. Also, since $H$ is sincere and $M(\alpha)$ 
is surjective, we have $M(a) \simeq H(a) \simeq k$. So $H(\alpha)$ is an isomorphism, a contradiction 
to its non-surjectivity. Since dim$_k M(s_i) = 1$, there exists an arrow $\beta \colon c \lra d$ on the 
path $(*)$ above such that 1 = dim$_k M(c) <$ dim$_k M(d)$. Therefore Coker$M(\beta) \neq 0$, hence 
$S_d$ is a summand of top$M$, a contradiction to top$M = S_{s_1} \oplus S_{s_2}$ and $ d \notin \{ s_1, s_2 \}$.
\end{proof}

\subsection{} 
{\sc Lemma. } {\it Let $C$ be a tame concealed algebra of type ${\widetilde{\D}_4}$
and $(H_{\lambda})_{\lambda}$ be an infinite family of pairwise non-isomorphic simple homogeneous modules. 
If $M$ is a module such that $H_{\lambda} \subset M$ for all $\lambda$, then $C[M]$ is wild. } 
\begin{proof} 
Assume first $C$ to be non-schurian, then $C$ is given by the quiver \\ 
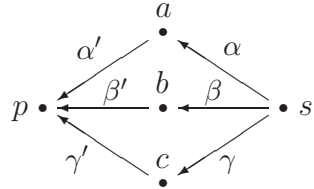
\begin{picture}(60,15)
\multiput(56,0)(16,0){3}{\circle*{1}}
\multiput(72,10)(0,-20){2}{\circle*{1}}

\multiput(70,0)(16,0){2}{\vector(-1,0){12}}
\multiput(70,9)(16,-10){2}{\vector(-3,-2){12}}
\multiput(70,-9)(16,10){2}{\vector(-3,2){12}}

\put(52,-1){$p$} \put(90,-1){$s$}
\put(71,12){$a$} \put(71,2){$b$} \put(71,-8){$c$} 
\put(80,7){\small$\alpha$}   \put(77.5,1){\small$\beta$}   \put(79.5,-8){\small$\gamma$} 
\put(60.5,6.5){\small$\alpha '$}   \put(64,1){\small$\beta '$}   \put(59,-8){\small$\gamma '$} 
\end{picture}
\vspace{1.4 cm}\\
bound by $\alpha \alpha' + \beta \beta' + \gamma \gamma' = 0$. For $\lambda \in k$, $H_\lambda$ is 
given by \\
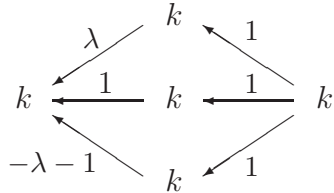
\begin{picture}(60,17)
\multiput(51,-1)(20,0){3}{$k$}
\multiput(71,10)(0,-22){2}{$k$}

\multiput(68,0)(20,0){2}{\vector(-1,0){12}}
\multiput(68,10)(20,-12){2}{\vector(-3,-2){12}}
\multiput(68,-10)(20,12){2}{\vector(-3,2){12}}

\put(81.5,8){\small 1}   \put(81.5,1){\small 1}   \put(81.5,-10){\small 1} 
\put(60,7){\small$\lambda$}   \put(62,1){1}   \put(50,-9.5){\small$-\lambda -1$} 
\end{picture}
\vspace{1.6 cm}\\
Since $H_\lambda \subset M$ for any $\lambda$, then $M$ is preinjective and we may further assume 
that dim$_k($top$M) \leq 2$. Also, soc$H_\lambda = S_p$ is a summand of soc$M$ while $S_s$ is a 
summand of top$M$. A straightforward examination of the preinjective component of $C$ shows that 
$M$ is indecomposable and is one of the four modules $I_p, \tau_C^2 I_a, \tau_C^2 I_b$ or 
$\tau_C^2 I_c$. In the first case, $C[I_p]$ has quiver  \\
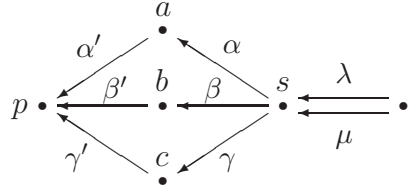
\begin{picture}(60,15)
\multiput(48,0)(16,0){4}{\circle*{1}}
\multiput(64,10)(0,-20){2}{\circle*{1}}

\multiput(62,0)(16,0){2}{\vector(-1,0){12}}
\multiput(94,1)(0,-2){2}{\vector(-1,0){12}}
\multiput(62,9)(16,-10){2}{\vector(-3,-2){12}}
\multiput(62,-9)(16,10){2}{\vector(-3,2){12}}

\put(44,-1){$p$} \put(79,2){$s$}
\put(63,12){$a$} \put(63,2){$b$} \put(63,-8){$c$} 
\put(72,7){\small$\alpha$}   \put(69.5,1){\small$\beta$}   \put(71.5,-8){\small$\gamma$} 
\put(52.5,6.5){\small$\alpha '$}   \put(56,1){\small$\beta '$}   \put(51,-8){\small$\gamma '$} 
\put(87,3){\small$\lambda$}  \put(87,-5){\small$\mu$} 
\end{picture}
\vspace{1.4 cm}\\
bound by $\lambda \alpha = \mu \alpha$, $\lambda \beta = \mu \beta$, $\lambda \gamma = \mu \gamma$, 
$\alpha \alpha' + \beta \beta' + \gamma \gamma' = 0$. Changing the presentation (replacing $\lambda$ 
by $\lambda - \mu$) gives the same quiver bound by $\lambda \alpha = 0, \lambda \beta = 0, \lambda \gamma = 0, 
\alpha \alpha' + \beta \beta' + \gamma \gamma' = 0$. This is a split extension \cite{ACT} of the 
algebra given by the quiver \\
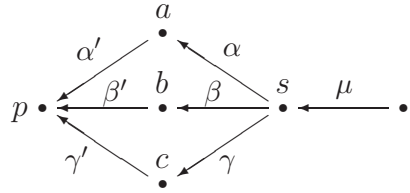
\begin{picture}(60,17)
\multiput(48,0)(16,0){4}{\circle*{1}}
\multiput(64,10)(0,-20){2}{\circle*{1}}

\multiput(62,0)(16,0){2}{\vector(-1,0){12}}
\multiput(94,0)(0,-2){1}{\vector(-1,0){12}}
\multiput(62,9)(16,-10){2}{\vector(-3,-2){12}}
\multiput(62,-9)(16,10){2}{\vector(-3,2){12}}

\put(44,-1){$p$} \put(79,2){$s$}
\put(63,12){$a$} \put(63,2){$b$} \put(63,-8){$c$} 
\put(72,7){\small$\alpha$}   \put(69.5,1){\small$\beta$}   \put(71.5,-8){\small$\gamma$} 
\put(52.5,6.5){\small$\alpha '$}   \put(56,1){\small$\beta '$}   \put(51,-8){\small$\gamma '$} 
\put(87,2){\small$\mu$}  
\end{picture}
\vspace{1.4 cm}\\
bound by $\alpha \alpha' + \beta \beta' + \gamma \gamma' = 0$, which is evidently wild. 

In the second case, $C[\tau_C^2 I_a]$ has quiver \\
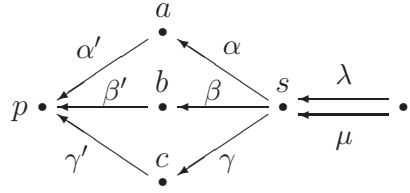
\begin{picture}(60,15)
\multiput(48,0)(16,0){4}{\circle*{1}}
\multiput(64,10)(0,-20){2}{\circle*{1}}

\multiput(62,0)(16,0){2}{\vector(-1,0){12}}
\multiput(94,1)(0,-2){2}{\vector(-1,0){12}}
\multiput(62,9)(16,-10){2}{\vector(-3,-2){12}}
\multiput(62,-9)(16,10){2}{\vector(-3,2){12}}

\put(44,-1){$p$} \put(79,2){$s$}
\put(63,12){$a$} \put(63,2){$b$} \put(63,-8){$c$} 
\put(72,7){\small$\alpha$}   \put(69.5,1){\small$\beta$}   \put(71.5,-8){\small$\gamma$} 
\put(52.5,6.5){\small$\alpha '$}   \put(56,1){\small$\beta '$}   \put(51,-8){\small$\gamma '$} 
\put(87,3){\small$\lambda$}  \put(87,-5){\small$\mu$} 
\end{picture}
\vspace{1.4 cm}\\
bound by $\lambda \beta = \mu \beta$, $\lambda \gamma = \mu \gamma$, 
$\alpha \alpha' + \beta \beta' + \gamma \gamma' = 0$. It contains as full convex subcategory the wild 
hereditary algebra with quiver \\
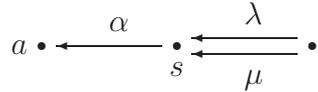
\begin{picture}(60,10)

\multiput(70,0)(0,-2){1}{\vector(-1,0){14}}
\multiput(88,1)(0,-2){2}{\vector(-1,0){14}}

\multiput(54,0)(18,0){3}{\circle*{1}}
\put(63,2){$\alpha$} 
\put(81,3){$\lambda$}  \put(81,-5){$\mu$}

\put(50,-1){$a$}  \put(71,-4){$s$}
\end{picture}
\vspace{1 cm}\\
The remaining cases follow by symmetry. 

We may thus assume that $C$ is schurian. If $C[M]$ is not wild, then $C$ is hereditary and has one of the 
quivers\\
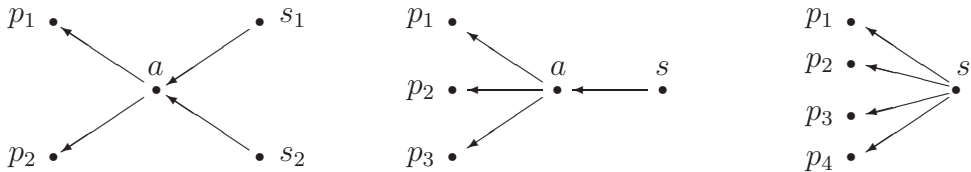
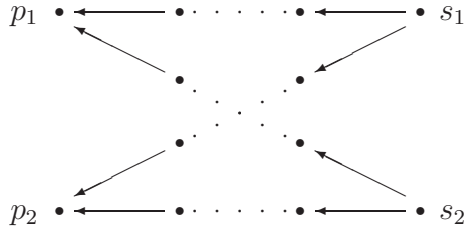
\begin{picture}(60,17)
\multiput(12,9)(27.4,0){2}{\circle*{1}}
\multiput(12,-9)(27.4,0){2}{\circle*{1}}
\put(25.7,0){\circle*{1}}
\multiput(38,8)(-14,-9){2}{\vector(-3,-2){11}}
\multiput(38,-8)(-14,9){2}{\vector(-3,2){11}}
\put(6,9){$p_1$}   \put(6,-9.5){$p_2$}  \put(42,9){$s_1$}   \put(42,-9.5){$s_2$}  \put(24.5,2){$a$}

\multiput(65,0)(14,0){3}{\circle*{1}}
\multiput(65,9)(0,-18){2}{\circle*{1}}
\multiput(77,0)(14,0){2}{\vector(-1,0){10}}
\put(77,1){\vector(-3,2){10}}   \put(77,-1){\vector(-3,-2){10}}
\put(59,9){$p_1$}   \put(59,-0.5){$p_2$}  \put(59,-9.5){$p_3$}   \put(92,2){$s$}  \put(78,2){$a$}

\multiput(118,9)(0,-18){2}{\circle*{1}} \multiput(118,3.5)(0,-7){2}{\circle*{1}} \put(132,0){\circle*{1}}
\put(131,1){\vector(-3,2){11}}  \put(131,0.5){\vector(-4,1){11}} 
\put(131,-0.5){\vector(-4,-1){11}}  \put(131,-1){\vector(-3,-2){11}} 
\put(112,9){$p_1$}   \put(112,3){$p_2$}  \put(112,-4){$p_3$}   \put(112,-10){$p_4$}  \put(132,2){$s$}
\end{picture}
\vspace{1.6 cm}\\
Recall that soc$H_\lambda$ is a summand of soc$M$ and dim$_k($top$M) \leq 2$. It is easily seen that 
in the first and third cases, no preinjective $C$-module satisfies these conditions. In the second case, 
however, we have the possibility $M = \tau^2_C I_s$, but then $C[M]$ is given by the quiver\\
\begin{picture}(60,17)
\multiput(48,0)(16,0){4}{\circle*{1}}
\multiput(48,10)(0,-20){2}{\circle*{1}}

\multiput(62,0)(16,0){2}{\vector(-1,0){12}}
\multiput(94,1)(0,-2){2}{\vector(-1,0){12}}
\put(62,1){\vector(-3,2){12}}
\put(62,-1){\vector(-3,-2){12}}

\put(42,10){$p_1$}    \put(42,0){$p_2$}     \put(42,-10){$p_3$}
\put(79,2){$s$}       \put(63,2){$a$}
\put(71,1){\small$\alpha$}  
\put(56,7){\small$\beta$}   \put(53.5,1){\small$\gamma$}   \put(55.5,-8.5){\small$\delta$} 
\put(87,3){\small$\lambda$}  \put(87,-5){\small$\mu$} 
\end{picture}
\vspace{1.6 cm}\\
bound by $\lambda \alpha \beta = \mu \alpha \beta$, $\lambda \alpha \gamma = \mu \alpha \gamma, 
\lambda \alpha \delta = \mu \alpha \delta$. It contains the wild hereditary full subcategory \\
\begin{picture}(60,10)

\multiput(70,0)(0,-2){1}{\vector(-1,0){14}}
\multiput(88,1)(0,-2){2}{\vector(-1,0){14}}

\multiput(54,0)(18,0){3}{\circle*{1}}
\put(63,2){$\alpha$} 
\put(81,3){$\lambda$}  \put(81,-5){$\mu$}

\put(50,-1){$a$}  \put(71,-4){$s$}
\end{picture}
\vspace{1 cm}\\
This establishes the lemma. 
\end{proof}

\pn {\sc Example.} Before proving our next result, we observe that the above lemma does not hold true 
for tame concealed algebras of type ${\widetilde{\A}}$. Indeed, let $C$ to be the Kronecker algebra, that is, the 
hereditary algebra given by the quiver \\
\begin{picture}(60,7)
\multiput(79,1)(0,-2){2}{\vector(-1,0){14}}

\multiput(63,0)(18,0){2}{\circle*{1}}

\put(59,-1){$p$}  \put(83,-1){$s$}
\end{picture}
\vspace{.7 cm}\\
Observe that the indecomposable injective $C$-module $I_p$ at the point $p$ contains all the 
indecomposable regular modules $H_\lambda$ given by \\
\begin{picture}(60,7)
\multiput(79,1)(0,-2){2}{\vector(-1,0){14}}

\multiput(63,0)(18,0){2}{\circle*{1}}

\put(59,-1){$k$}  \put(83,-1){$k$}
\put(72,3){$1$}  \put(72,-5){$\lambda$}
\end{picture}
\vspace{.7 cm}\\
Now, the one point extension $C[I_p]$ is the radical square zero algebra given by the 
quiver \\
\begin{picture}(60,7)
\multiput(70,1)(0,-2){2}{\vector(-1,0){14}}
\multiput(88,1)(0,-2){2}{\vector(-1,0){14}}

\multiput(54,0)(18,0){3}{\circle*{1}}
\put(63,3){$\beta$}  \put(63,-5){$\delta$}
\put(81,3){$\alpha$}  \put(81,-5){$\gamma$}
\put(50,-1){$p$}  \put(71,-4){$s$}
\end{picture}
\vspace{.7 cm}\\
bounded by $\alpha \delta = 0$, $ \gamma \beta = 0 $ and $\alpha \beta = \gamma \delta$, 
which is clearly tame. 
\vs

\pn {\sc Remark.} Notice, however, that if $C$ is tame concealed of type ${\widetilde{\A}}$ and schurian 
and $M$ satisfies the conditions of (2.3), then $C[M]$ is wild. Indeed, if this is the case then, 
because of (2.2), $C$ has at most 2 sources, hence is hereditary with quiver \\
\begin{picture}(60,25)
\multiput(48,12.5)(16,0){4}{\circle*{1}}
\multiput(64,3.5)(16,0){2}{\circle*{1}}
\multiput(48,-14)(16,0){4}{\circle*{1}}
\multiput(64,-5)(16,0){2}{\circle*{1}}

\multiput(62,12.5)(32,0){2}{\vector(-1,0){12}}
\put(62,4.5){\vector(-2,1){12}}
\put(94,10.5){\vector(-2,-1){12}}
\multiput(62,-14)(32,0){2}{\vector(-1,0){12}}
\put(62,-6){\vector(-2,-1){12}}
\put(94,-12){\vector(-2,1){12}}

\multiput(66,12.5)(3,0){5}{\circle*{.1}}
\multiput(66,-14)(3,0){5}{\circle*{.1}}
\multiput(66,-4)(3,1.5){5}{\circle*{.1}}
\multiput(78,-4)(-3,1.5){5}{\circle*{.1}}

\put(41.5,11.5){$p_1$}
\put(98.5,11.5){$s_1$}
\put(41.5,-15){$p_2$}
\put(98.5,-15){$s_2$}
\end{picture}
\vspace{2.5 cm}\\
Here, soc$H_\lambda = S_{p_1} \oplus S_{p_2}$ is a summand of soc$M$ while top$M = $ top$H_\lambda = S_{s_1}\oplus S_{s_2}$. It is easily seen that no preinjective $C$-module satisfies these conditions. 

\subsection{} 
{\sc Lemma.} {\it Let $C$ be a tame concealed algebra of type ${\widetilde{\D}_n}$, with $n \geq 5$, 
or ${\widetilde{\E}}$ and $H$ be a simple homogeneous $C$-module. If $M$ is preinjective such that 
$H \subset M$, then $C[M]$ is wild. } 
\begin{proof} 
We claim that $C[M]$ is wild. Indeed, the tubular type of $C$ is of the form $(a,b,c)$ with at 
least one of $a, b, c$ larger than or equal to 3. Taking then a one-point extension by a simple 
homogeneous module yields a wild algebra. This establishes our claim. Applying now (1.2), we get
that $C[M]$ is wild, as required. 
\end{proof}

\subsection{} 
We are now able to prove the main result of this section. 
\vs\\
{\sc Theorem.} {\it Let $C$ be a tame concealed algebra of type distinct of ${\widetilde{\A}}$, and 
$(H_\lambda)_\lambda$ be an infinite 
family of pairwise non-isomorphic modules. Suppose $M$ is a $C$-module such that $H_\lambda \subset M$ 
for all $\lambda$. Then $C[M]$ is wild. } 
\begin{proof} This follows immediately from (2.3) and (2.4). 
\end{proof} 

\section{Torsionless finiteness and representation dimension.} 

\subsection{} 
Let $A$ be an Artin algebra. An $A$-module $M$ is called a {\it generator} (of mod$A$) if $A_A \in$
add$M$ and a {\it cogenerator} if (D$A)_A = $ Hom$_k(_AA,k) \in $ add$M$. Let $A$ be a non-semisimple 
algebra. The {\it representation dimension} rep.dim$.A$ of $A$ is the infimum of the global dimensions 
of the algebras End$_AM$, where $M$ is a generator and a cogenerator of mod$A$, see \cite{Au}. 

An Artin algebra $A$ is called {\it torsionless-finite} if every indecomposable projective $A$-module
has only finitely many isomorphism classes of indecomposable submodules. One defines {\it cotorsionless-finite}
dually.  An Artin algebra is torsionless-finite if and only if it is cotorsionless-finite \cite{AB}. 
We need essentially the following result (see \cite{Ri3}). 
\vsd\\
{\sc Theorem.} {\it If $A$ is torsionless-finite, then rep.dim.$A \leq 3$. }

\subsection{} 
{\sc Lemma.} {\it Let $A$ be a (possibly wild) branch enlargement of a tame concealed algebra. Then
$A$ is torsionless-finite. In particular, rep.dim.$A \leq 3$. } 
\begin{proof} 
Using the description of mod$A$ in \cite{LS}, we see that, if $P_A$ is an indecomposable projective 
$A$-module, then either $P_A$ is postprojective (in which case it clearly has only finitely many 
isomorphism classes of indecomposable modules) or it lies in an inserted tube $\Gamma$. But in this 
latter case, $P$ has only finitely many indecomposable submodules lying in $\Gamma$ and, using \cite{AR}, 
there are only finitely many postprojective modules $X$ such that dim$_k X \leq $ dim$_k P$. Since 
the other tubes in the same family as $\Gamma$ are orthogonal to $\Gamma$, the proof is complete. 
\end{proof} 

\subsection{} 
{\sc Corollary.} {\it Let $A$ be a tame quasi-tilted algebra, then $A$ is torsionless-finite. 
In particular, rep.dim.$A \leq 3$. } 
\begin{proof}
This follows from (3.2) and \cite{Sk3} 
\end{proof} 

\vsd\pn
{\sc Remark.} This corollary is a particular case of the main result of \cite{Op}. 
Further, the same proof as in (3.2) and \cite{AS1} give that if $A$ is iterated tilted of euclidean 
type, then $A$ is torsionless-finite, and so rep.dim.$A \leq 3$. This is a particular case of the main 
result of \cite{CHU}.

\subsection{} 
{\sc Lemma.} {\it Let $A$ be a strongly simply connected tame coil enlargement of a tame concealed algebra, 
then $A$ is torsionless-finite. In particular, rep.dim.$\leq 3$. } 
\begin{proof} 
By (1.3), $A$ contains as full convex subcategories a unique maximal branch coextension $A^-$ and a 
unique maximal branch extension $A^+$ of a tame concealed algebra, and mod$A$ contains a weakly separating 
family of coils ${\mathcal T}'$ such that ind$A = {\mathcal P}' \vee {\mathcal T}'\vee {\mathcal Q}'$ 
where ${\mathcal P}' \subset $ ind$A^-$ and ${\mathcal Q}' \subset $ ind$A^+$. Let $P_A$ be an 
indecomposable projective module. We have three cases:

(a) Assume $P \in {\mathcal P}'$, then $P$ is an indecomposable projective $A^-$-module. Because of 
(3.2), and the fact that $A^-$ is a branch enlargement of a tame concealed algebra, $P$ has only 
finitely many isomorphism classes of indecomposable submodules in mod$A^-$. Now, the indecomposable 
submodules of $P$ in mod$A$ and mod$A^-$ coincide. 

(b) Assume $P \in {\mathcal Q}'$. For the same reason, $P$ has only finitely many isomorphism classes 
of indecomposable submodules in mod$A^+$, hence in mod$A$. 

(c) Assume $P \in {\mathcal T}'$, and let $M$= rad$P$. There exists a sequence of full convex 
subcategories of $A$
$$ A^- \ = \ A_0 \ \subset\ A_1 \ \subset \ \cdots \ \subset \ A_t = A$$
which are iterated one-point extensions of $A^-$, and an index $i$ such that $A_{i+1} = A_i[M]$, that is, 
$P$ is the unique indecomposable projective in mod$A_{i+1}$ which is not in mod$A_i$. Also, $M$ is the 
pivot of an admissible operation and then is indecomposable except perhaps in the case (ad1). Then, 
$M = M' \oplus M''$, where $M''$ is a directed indecomposable, while $M'$ is an indecomposable lying in  
a coil. Since $M''$ has only finitely many isomorphism classes of indecomposable submodules and $A_{i+1}$
is of the form $(A_i'\times D)[M'\oplus M'']$, where $D$ is a triangular matrix algebra, and $A_i'$ is the 
full subcategory of $A_i$ with objects $(A_i')_0 = (A_i)_0 \setminus D_0$, then it suffices to consider the 
submodules of $M'$. We may thus for simplicity assume that $M = M'$. \\
Let $M^-$ be the largest $A^-$-submodule of $M$. Assume that $M^-$ has infinitely many isomorphim classes 
of indecomposable $A^-$-submodules. Because of (3.2), $M^-$ has to be a preinjective $A^-$-module. By (2.5), 
$A^-[M^-]$ is wild. Since $A^-$ is a full convex subcategory of $A_i$, then $A_i[M^-]$ is wild. By (1.2), 
and the tameness of $A_i$, we get that $A_{i+1} = A_i[M]$ is wild. This is absurd, because $A_{i+1}$ is a 
full convex subcategory of $A$ which is tame. This shows that $M^-$ has only finitely many non-isomorphic 
indecomposable submodules. Since a submodule of $M$ is either a submodule of $M^-$, or lies in the coil 
containing $M$, we are done. 
\end{proof} 

\subsection{} 
{\sc Corollary.} {\it Let $A$ be a multicoil algebra, and $P$ be an indecomposable projective $A$-module 
lying in a coil $\Gamma$ of $\Gamma($mod$A$), then $P$ has only finitely many isomorphism classes of 
indecomposable submodules. } 
\begin{proof} 
By \cite{AS3}(5.9), we can assume that $\Gamma$ is obtained from a stable tube over a tame concealed algebra
$C$ by a sequence of admissible operations and the support algebra $B$ of $\Gamma$ is obtained from $C$ by the
corresponding sequence of one-point extensions and coextensions. Since the $A$-submodules and the $B$-submodules
of $P$ coincide, the result follows from (3.4). 
\end{proof} 

\subsection{} 
We are now able to state and prove the main result of this paper. 
\vsd\\
{\sc Theorem.} {\it Let $A$ be a strongly simply connected algebra of polynomial growth, then 
$A$ is torsionless-finite. In particular, rep.dim.$A \leq 3$.} 
\begin{proof} 
We may clearly assume that $A$ is representation-infinite. Let $P$ be an indecomposable projective 
$A$-module. By (1.6), if $P$ is not directed, then the support algebra $B$ of $P$ is a tame coil 
enlargement of a tame concealed algebra. By (3.4), $P$ has only finitely many isomorphism classes of 
indecomposable $B$-submodules, hence $A$-submodules. On the other hand, if $P$ is directed, then, by 
(1.6) again, $B$ is tame and tilted. Since, clearly, $B$ has a sincere directed indecomposable module 
(namely $P$), then, by \cite{Pe}, $B$ is domestic in at most two one-parameters. Moreover, by \cite{Bo}, 
$B$ is a full convex subcategory of $A$. \\
Let $\Gamma$ denote the component of $\Gamma$(mod$B$) containing $P$, and $\Sigma$ denote the full 
subquiver of $\Gamma$ consisting of the indecomposable modules $X \in \Gamma$ such that there is a 
path $X \leadsto P$, and every such path is sectional. By \cite{ASS}(IX.2.6, p. 364), $\Sigma$ is a 
complete slice in $\Gamma$(mod$B$). Note also that in $\Gamma($mod$B$), there are finitely many directed 
components preceding $\Gamma$ (actually, by \cite{Pe}, at most two). Since indecomposable modules lying 
in a directed component are uniquely determined by their composition factors \cite{AR}, then $P$ has 
infinitely many isomorphism classes of indecomposable submodules if and only if there exists an 
infinite family $({\mathcal T}_\lambda)_\lambda$ of homogeneous tubes over a tame concealed algebra $C$, 
and an infinite family of homogeneous modules $H_\lambda \in {\mathcal T}_\lambda$ such that 
$H_\lambda \subset $ rad$P$, for each $\lambda$. Clearly, then, all the $H_\lambda$ are contained in the 
largest $C$-submodule $R$ of rad$P$. Moreover, by the structure of homogeneous tubes, we can assume 
the $H_\lambda$ to be simple homogeneous. By (2.5), if this is the case, then $C[R]$ is wild, hence 
so is $B$. This shows that $P$ has only finitely many isomorphism classes of indecomposable $B$-submodules, 
hence $A$-submodules.
\end{proof}

\subsection{} 
{\sc Remark.} As is seen from the example in Section 2, it is not true in general that multicoil algebras 
are torsionless-finite. However, notice that the algebra \\
\begin{picture}(60,7)
\multiput(70,1)(0,-2){2}{\vector(-1,0){14}}
\multiput(88,1)(0,-2){2}{\vector(-1,0){14}}

\multiput(54,0)(18,0){3}{\circle*{1}}
\put(63,3){$\beta$}  \put(63,-5){$\delta$}
\put(81,3){$\alpha$}  \put(81,-5){$\gamma$}
\end{picture}
\vspace{.7 cm}\\
bounded by $\alpha \delta = 0$, $ \gamma \beta = 0 $ and $\alpha \beta = \gamma \delta$ 
is tilted and so, by \cite{APT}, has representation dimension at most 3.

\makeatletter
\renewcommand{\@biblabel}[1]{\hfill#1.}\makeatother

\end{document}